\documentclass{article}
\usepackage{cite}
\title{On f-Derangements and Decomposing Bipartite Graphs into Paths}
\author{Michael Plantholt\\
Department of Mathematics\\
Illinois State University\\
Normal, IL 61790-4520, USA\\
Tel: 309-261-2039\\
email: mikep@ilstu.edu\\
\and
Hamidreza Habibi\\
Department of Economics\\
University of California, Santa Cruz\\
Santa Cruz, CA 95064, USA 
\and
Benjamin Mussell \\
Department of Mathematics\\
Illinois State University\\
Normal, IL 61790-4520, USA\\}

\date{\today}

\let\oldbibliography\thebibliography
\renewcommand{\thebibliography}[1]{%
  \oldbibliography{#1}%
  \setlength{\itemsep}{-2pt}%
}

\begin{document}
\maketitle
\newtheorem{theorem}{Theorem}
\newtheorem{lemma}{Lemma}
\newenvironment{proof}{{\bf Proof}.}{\hspace{3mm} \rule{3mm}{3mm}}
\newtheorem{conjecture}{Conjecture}

\begin{abstract}

Let $f: \{1, ..., n\} \rightarrow \{1, ..., n\}$ be a function (not necessarily one-to-one).  An $f-derangement$ is a permutation  $ g:\{1,...,n\} \rightarrow \{1,...,n\}$ such that $g(i) \neq f(i)$ for each  $ i = 1, ..., n$.  When $f$  is itself a permutation, this is a standard derangement.  We examine properties of f-derangements, and show that when we fix the maximum number of preimages for any item under $f$, the fraction of permutations that are f-derangements tends to $ 1/e$ for large $n$, regardless of the choice of $f$.  We then use this result to analyze a heuristic method to decompose bipartite graphs into paths of length 5.\\

  \medskip
keywords: derangement, decomposition, matching, directed path

\end{abstract}

\section{Introduction}

Recall that a $derangement$ is a permutation  $ g:\{1,...,n\} \rightarrow \{1,...,n\}$ such that $g(i) \neq i$ for each $ i = 1, ..., n$. The classic problem of enumerating derangements was first formulated in 1708 by P. R. de Montmort, and solved by him in 1713.  At about the same time Nicholas Bernoulli also solved the problem using inclusion-exclusion.  A classic result(see, for example, \cite{Wells}) is that as $n$ tends to infiniy, the fraction of permutations that are derangements tends to $ \frac{1}{e} $. Now let  $f: \{1, ..., n\} \rightarrow \{1, ..., n\}$ be a function (not necessarily one-to-one).  An $f-derangement$ is a permutation  $ g:\{1,...,n\} \rightarrow \{1,...,n\}$ such that $g(i) \neq f(i)$ for each  $ i = 1, ..., n$.  It is clear that when $f$  is itself the identity function, or more generally any permutation, an f-derangement is equivalent to a standard derangement.  If function $f$ has at most $k$ items mapped to the same value, we will call $f$ a $k-max$  function.

In this section, we provide some elementary results on f-derangements.  In Section 2 we focus on  f-derangements of a 2-max function $f$, giving some recursive formulas that can be used to find the number of them, and use these to show that for large $n$ and any 2-max f-function, the percent of permutations on $n$ elements that are f-derangements is approximately $\frac{1}{e}$.  In Section 3, we extend that result to k-max f-derangements, showing the following.\\ \medskip

Theorem 1.  Let $k$ be a fixed positive integer, and let $f_{1}, f_{2}, f_{3}, ... $ be a sequence of k-max functions where $f_{i}: \{1, ..., i\} \rightarrow \{1, ..., i\}$.  Then as $n \rightarrow \infty$, the percentage of permutations on $n$ elements that are $f_{n}$-derangements has limit $ \frac{1}{e} $.  \\

Finally, in Section 4 we present a heuristic algorithm that seeks to partition a 5-regular bipartite graph into paths of length 5.  We use the Theorem above to get a measure for how frequently we would expect that algorithm to succeed.

For basic properties of f-derangements, we note that for a given function  f, there may be no f-derangements; indeed note that if $f(i) = 1$ for all values of $i$, no f-derangement exists. For a less restrictive example, suppose that $ f(1) = 2$ and$ f(i) = 1$ for all other values of $i$. Then any f-derangement g must have g(1) = 1, and all other values are then unrestricted, so there are $ (n-1)! $ such f-derangements. 

Next suppose that $n$ is even, and $f(i)=1$ for $ i = 1, 2, ..., n/2$, while $ f(i) = 2 $ otherwise. It is easy to see then the the number of f-derangements will be exactly $ \frac{n}{2} \cdot \frac{n}{2} \cdot (n-2)! $.  Thus the percent of permutations that will be f-derangements for such a function will be 
$\frac{ \frac{n}{2} \cdot \frac{n}{2}}{(n \cdot (n-1))} $, which approaches 1/4 as $n$ gets large.  It is not hard to establish also that when we split the images of our function f equitably among 3 options, the limiting value is $(2/3)^{3}$, and more generally an equitable split between  s  images will have $ ( \frac{(s-1)}{s})^{s}$  as the limiting value for the percent of permutations that are f-derangements.  Note that for $s$ large this will be approximately $1/e$, the famous limiting value for derangements.  This bound will appear again when we turn to k-max derangements.

\section{2-max f-derangements}

We begin by presenting tables for the number of 2-max f-derangements for different 2-max functions for permutations of size 4, 5, and 6.  To simplify, we note that the numbers of such f-derangements do not depend on specific values mapped to, but only the number of collision pairs in the mapping.  For example, both the 2-max functions $f_1$ and $f_2$ will have the same number of f-derangements if  $f_1(1)= 1, f_1(2) = 2, f_1(3)=f_1(4) = 3$ and  $f_2(1) = 1, f_2(2) = f_2(3) = 2, f_2(4) = 4$  . So we let $A$ denote the number of items that have a collision in their mapping (so $A$ is always even in a 2-max setting), and $B$ denote the number of items that have a unique image.  Thus both $f_1$ and $f_2$ above have values $A=2, B=2$.  Finally, for this construction we will extend our requirements on f to allow a function that is undefined on some values of $ \{ 1, 2, ...n \}$; we will let  $C$ denote that number in the 2-max case.  Thus the function $f$ on $n=6$ items with $ f(1) = 1, f(2)=2, f(3)=2, f(4)= 4, f(5)=5, f(6)$ undefined has parameters $ A=2, B=3, C=1$.  Tables 1, 2 and 3 below give the number of 2-max f-derangements for $ n= 4, 5, 6 $.  We will derive these values and properties of them in what follows.

\begin{table}[t]
\caption{number of 2-max derangements for $n=4$ }
\centering
\begin{tabular} {c c c c}
\hline
n=4 & A=0 & A=2 & A=4 \\[0.5ex]  
\hline
B=0 & 24 & 12 & 8  \\
B=1 & 18 &10  \\
B=2 & 14 & 8  \\
B=3 & 11  \\
B=4 & 9  \\ [1ex]
\hline
\end{tabular}
\label{table: f-derangements, n=4}
\end{table}

\begin{table}[h]
\caption{number of 2-max derangements for $n=5$ }
\centering
\begin{tabular} {c c c c}
\hline
n=5 & A=0 & A=2 & A=4 \\[0.5ex]  
\hline
B=0 & 120 & 72 & 48  \\
B=1 & 96 &60 & 40  \\
B=2 & 78 & 50  \\
B=3 & 64 & 42  \\
B=4 & 53  \\
B=5 & 44 \\ [1ex]
\hline
\end{tabular}
\label{table: f-derangements, n=5}
\end{table}

\begin{table}[h]
\caption{number of 2-max derangements for $n=6$ }
\centering
\begin{tabular} {c c c c  c}
\hline
n=6 & A=0 & A=2 & A=4 & A=6 \\[0.5ex]  
\hline
B=0 & 720 & 480 & 336 & 240  \\
B=1 & 600 & 408 & 288  \\
B=2 & 504 & 348 & 248  \\
B=3 & 426 & 298  \\
B=4 & 362 & 256  \\
B=5 & 309 \\
B=6 & 265 \\ [1ex]
\hline
\end{tabular}
\label{table: f-derangements, n=6}
\end{table}

Such tables for 2-max derangements are easily constructed with the help of recursive equations.  We provide two such equations.  Let  $f$  be a 2-max function from some subset of $\{ 1, 2, ..., n\}$ to $\{ 1, 2, ..., n\}$ .  Using the notation given above, if $A =$ the number of items which map to the same elements as some other item, $B$ = the number of items that has only one pre-image, and $C = n - A - B$ is the number of items that are not mapped by $f$, we dentote the number of f-derangements by D[A,B,C], and sometimes denote $f$ itself by f[A,B,C]. \\

Recursion Formula 1:   If $ A, B > 0$ then \\
$ D[A,B,C] =( C+ \frac{A}{2}) \cdot D[A,B-1,C] + (B-1) \cdot D[A, B-2, C+1]+ \frac{A}{2}] \cdot D[A-2, B-1, C+2]$ \\
\\

The derivation of Recursion Formula 1 is straightforward.  Without loss of generality, assume that $f(1) = 1$ and no other element maps to 1 under $f$. For a permutation $g$ to be an f-derangement, we must have that $g(1) \neq 1$.  So $g(1)$ could be any of the $C + \frac{A}{2}$ elements that are not mapped to by any element of $f$, and of the other $B-1$ elements that are mapped to by a single element in $f$, or any of the $ \frac{A}{2} $ elements that have a preimage of cardinality 2 under $f$.  This breakdown gives the first item in the each of the three terms in the recursion sum.  Now suppose that $g(1) = i$, where $i$ is not the image of any element of $f$.  To complete the options for $g$ to be an f-derangement, we now need to assign the remaining $n-1$ elements.  We still have $A/2$ pairs that need to avoid the same element, now have $B-1$ elements that have a unique term to avoid in the mapping, and $C - 1 + 1$ elements that no item needs to avoid.  Thus the first term in the sum is  $(C+ \frac{A}{2}) \cdot D[A,B-1,C]$ ; the other terms follow by similar reasoning. \\

Thus, using our $n=5$  tables to get $D[2,3,1]$ with $n=6$, we have $D[2,3,1] = 2 \cdot D[2,2,1] + 2 \cdot D[2, 1, 2] + 1 \cdot D[0, 2, 3] = 2 \cdot 50 + 2\cdot 60 + 1 \cdot 78 = 298 $  \\

Our next recursion formula is less obvious, but useful in determining the limiting number of f-derangements. \\

Recursion Formula 2:  If $B \geq 2$ then  $D[A,B,C] = D[A+2, B-2, C] + D[A,B-2, C]. $ \\
\\
To see why this formula holds, we begin by assuming $f$ is a 2-max function with the parameters $ A, B, C$.  We can assume without loss of generality that $f(1)=1$ and $f(2)=2 $, and no other elements are mapped to 1 or 2 by $f$.  The f-derangements can be partitioned into two types.  First, such a derangement $g$ could have $g(1)=2, g(2)=1$.  It is easy to see the number of these is $D[A, B-2, C]$.  The remaining possibilities for f-derangements have either one or both of $\{1,2\}$ mapping into $\{3, ..., n\}$.  Let $f^*$ be the function that is identical to $f$, except that $f^*(2) = 1$.  We claim that the number of  f-derangements that have either one or both of $\{1,2\}$ mapping into $\{3, ..., n\}$ is equal to the number of $f^*-derangements$. Certainly the number of $f, f^*$ derangments such that neither of the elements $\{1, 2\}$ maps into $\{1, 2\}$ is identical.  Finally, suppose we have an $f-derangement$, call it $g$  in which exactly one of $\{1, 2\}$ maps into $\{1, 2\}$.  If $g(1)=2, g(2) =i$, we get a corresponding $f^*-derangement $, call it $g*$, which is identical to $g$ except that $g^*(1)=i, g^*(2)=2$.  If $g(1)=i, g(2) =1$, we let $g^*=g$.  This correspondence shows equivalence between the number of $f$ and $f^*$ derangements in this partition, and thus that number is D[A+2,B-2,C].  The recurrence follows. \\
\\

As an example, using our tables for $n=4$ and $n=6$, we get that $D[2, 3, 1] = D[4, 1,1] + D[2,1,1] = 288 + 10 = 298$.

We note 3 properties from the tables, and give explanations for them. \\

Property 1.  For any 2-max function $f$ with parameter value $A > 0$, the number of f-derangements is even.
To see why this is so, without loss of generality suppose that $f(1) = f(2)$.  Then for any f-derangement $g$ for which $g(1) = i, g(2) = j $ we have a corresponding f-derangement $g^*$ obtained from $g$ by merely interchanging the images of 1 and 2.  Thus the f-derangements come in pairs, so the number of them is even.

Of course as we move down vertically or to the right horizontally through the tables, the number of f-derangements decrease, because we are adding more restrictions. But note that as we move diagonally from top right to top left, the values increase. We state this more formally as follows.\\
\\

Property 2.  For any 2-max function $f$ and parameter value $A > 0$,  $D[A, B, C] > D[A+2, B-2, C]$.
This follows immediately from Recursion Formula 2, because  $D[A,B,C] = D[A+2, B-2, C] + D[A,B-2, C]. $  \\

Finally recall the famous result that for large $n$, the fraction of permutations that are derangements is approximately $1/e$, or approximately $0.36$.  The number of regular derangements in the charts is given by $D[0,n,0]$.  In Table 3 with $n=6$, notice that at the other extreme, we have  $D[6,0,0]/n! = 240/720 = 1/3$, not far from $1/e$.  \\
\\
Property 3.  For $n$ large, for any 2-max $f: \{1, ..., n\} \rightarrow \{1, ..., n\}$, the fraction of permutations on $n$ items that are f-derangements is approximately $1/e$. \\
\\
To verify this, let $ D(n) $  denote the number of derangements of order $n$.  Let $f$ be a 2-max function with parameters  $f[A,B,0]$. If $A>0$, then from Recursion Formula 2, $D[A, B, 0] = D[A-2, B+2,0] - D[A, B-2,0]$.  Iterating again if $ A-2 > 0$ we get  $D[A, B, 0] = D[A-4, B+4] - D[A, B-2,0] - D[A-2, B, 0].$  After $A/2$ applications of the recursive formula, we get $D[A, B, 0] = D(n) $ minus $A/2$ terms, each at most $ (n-2)!.$  Thus as $n \rightarrow \infty , \frac{D(n) - D[A,B,0]}{n!} \rightarrow 0$.  Therefore as $ n \rightarrow \infty, \lim \frac{D(f)}{n!}  = \lim \frac{D(n)}{n!} = 1/e $.  \\
\\
In the next Section, we prove a theorem generalizing this result to include all k-max functions for k fixed.

\section{Limiting Value for k-max Derangements}

We seek to extend the previous result on the limiting behaviour of the number of derangements for 2-max functions to k-max functions for a fixed $k$.  In particular, we show the following.\\

Theorem 1.  Let $k$ be a fixed positive integer, and let $f_{1}, f_{2}, f_{3}, ... $ be a sequence of k-max functions where $f_{i}: \{1, ..., i\} \rightarrow \{1, ..., i\}$.  Then as $n \rightarrow \infty$, the percentage of permutations on $n$ elements that are $f_{n}$-derangements has limit $ \frac{1}{e} $.  \\

We first prove the following.

Lemma.  Let $f, f^*$ be functions with  $f, f ^* : \{1, ..., n\} \rightarrow \{1, ..., n\}$.  Suppose $f, f^*$ are identical except that $ f(i) = 1$ for $i = 1, 2, ..., r  (r \geq2)$ and no values are mapped to 2 under $f$, while $f^*(r)=2 $.  Then:\\
    (*)  $D[f] \leq D[f^*] $ , and\\
    (**)  $ D[f^*] - D[f] \leq  (n-2)! $\\
\\

Proof.  We focus on the possibilities of the images of elements $1, r$  in valid derangements  $g$ under functions $ f, f^*$ .\\
Case 1:  Neither of $1,r$ gets mapped into $ \{ 1,2 \}$ by $g$.
Because the only difference between $f,f^*$ is that $f(r) = 1, f^*(r)=2$, it is clear that the number of $f, f^* $ -derangements with this additional condition is identical.

Case 2:  Exactly one of $ \{ 1, r \}$ gets mapped into $ \{ 1,2 \}$ by $g$.
With this restriction, in order for $g$ to be an f-derangement, we must have $g(1) = 2, g(r) > 2$ or $g(r)=2, g(1) > 2  $. The number of  f-derangements of this type with $g(1) = 2, g(r) > 2$  is equal to the number of $ f^*$-derangements with the same $g-values$.  On the other hand, the number of  f-derangements of this type with  $g(r)=2, g(1) > 2  $  is the same as the number of $ f^*$-derangements with $ g(r) = 1, g(1) > 2 $.  Thus there will be the same number of $f, f^*$-derangements in this case.

Case 3: Both elements of $ \{ 1, r \}$ gets mapped into $ \{ 1,2 \}$ by $g$.
Such a permuation $g$ cannot be an f-derangment, because this forces either $g(1) = 1 $ or $ g(r) = 1$.  However, there may be $f^*$-derangements in this case that have $g(1)=2, g(r)=1 $.  Note however that the number of these that are $f^*$-derangements cannot be more than $ (n-2)! $.

The results now follow.

\section{Decomposition of Bipartite Graphs into Paths}

We now show how the previous results help in the analysis of a heuristic algorithm for decomposing bipartite graphs into paths. Recall that a bipartite graph is a graph in which there is a partition of the vertices into sets $ X, Y$  such that each edge of the graphs joins a vertex in  $X$  to a vertex in $Y$.  For additional graph theory terminoloy we follow \cite{West}.  One of the most famous unproven conjectures in graph theory is by Ringel \cite{Ringel}.\medskip

Conjecture (Ringel's Conjecture 1963).  For any tree  $T$  with  $s$ edges, the complete graph $K_{2s+1}$ has a decomposition into subgraphs, each isomorphic to $ T$ .\medskip

A bipartite version of this is the following, whose origin is unclear. \medskip

Conjecture (Folklore).  For any tree $T$ with $s$ edges, there is a decomposition of the complete bipartite graph $ K_{s,s} $  into subgraphs, each isomorphic to $T$.

Ringel's Conjecture was recently proved for sufficiently large graphs \cite{Mont} but both these Conjectures appear to be far from fully proved. There have been hundreds or even thousands of papers devoted to them.  Graham and Haggkvist \cite{Graham} have proposed a further strengthening of the conjectures, to include decompositions of general regular graphs (not necessarily complete).\medskip

Conjecture:   For any tree $T$ with $s$ edges, \\
	1. every 2s-regular graph has a decoposition into isomorphic copies of $T$. \\
	2.  every s-regular bipartite graph has a decomposition into isomorphic copies of $T$.\medskip

We focus on the problem of decomposing the edges of a bipartite $s$-regular graph into paths of length $s$.  That this works for any 3-regular bipartite graph follows readily from independent results by Bouchet \& Fouquet \cite{Bouchet}, and Kotzig \cite{Kotzig}.  The result was verified for 4-regular bipartite graphs by Jacobson, Truszczynski and Tuza in 1991 \cite{Jacobson}.  More recently, it was proved for 5-regular bipartite graphs (and more generally 5-regular triangle-free graphs) by Botler, Mota and Wakabayashi \cite{Botler}.  The proof by Botler et al. consists of starting with an approximate decomposition, then showing it is always possible to make adjustments to this approximation which will yield a valid decomposition.  We offer a different heuristic approach, which is simpler in concept, to try building such a decomposition, and analyze how frequently this will be successful.

We give our algorithm as follows.  Let  $G$ be a 5-regular bipartite graph, with partite sets  $X, Y$  such that each edge connects a vertex of  $X$  with a vertex of $Y$.  \medskip

Step 1:  Find a perfect matching $ M$ of $ G$, and let $G^*$ denote the 4-regular graph obtained by deleting the edges of $M$ from $G$.  Such a perfect matching exists by Konig's Theorem; in fact, we note for later discussion that by Konig's Theorem  $G$ has a decomposition into 5 perfect matchings.\medskip

Step 2:  Get a decomposition of $G^*$ into paths of length 4, which is possible by the result of Jacobson et al.  In fact, El-Zanati and Plantholt (private communication) showed that by tweaking the proof in that paper, we can guarantee that in the decompostion, we can direct the paths and guarantee that each vertex of $X$ is the start vertex of one path, the center vertex of one path, and the end vertex of one path. We assume our decomposition has that property.\medskip

Step 3:  To each path of our decomposition in Step 2, we attach to the start the edge of $M$ that includes the $X$-vertex that begins the path.\medskip

The three steps above give us a decomposition into subgraphs with 5 edges each.  If we are lucky, each such subgraph will be a path, and we are done.  For example, suppose our start graph is $K_{5,5}$, with X-vertices  $ 1,2,3,4,5$, and Y-vertices $A,B,C,D,E $, and let our perfect matching $M$ be the edges $ 1A, 2B, 3C, 4D, 5E$ .One possible decomposition of $G^*$ into paths of length 4 satisfying our desired conditions is given by the following 5 paths (note when strung together they give an euler tour):\\
$1,B,3,D,5$ \\
$5,C,1,D,2$\\
$2,A,4,E,3$\\
$3,A,5,B,4$\\
$4,C,2,E,1$\\
Attaching the edges of $M$ given by $A1$ to the start of the first of these, $E5$ to the second, $B2$ to the third,  $C3$ to the fourth and $D4$ to the last then gives us the desired decomposition of $K_{5,5}$ into paths of length 5.

Now, on the other hand, suppose our decomposition of $G^*$ is given by the following 5 paths:\\
$ 1,B,5,A,2$ \\
$2,D,3,A,4$ \\
$ 4,E,2,C,1$ \\
$3, E, 1, D, 5$ \\
$ 5, C, 4, B, 3$  \\

Now if we try extending these to paths of length 4 by adding $A1$ at the start of the first path, $B2$ to start the second, $D4$ the third, $C3$ the fourth and $E5$ the last, we see the first of these subgraphs is not a path, because we added $A1$ to a path that already had $A$ as its fourth vertex.  So the algorithm in this case does not yield a decomposition into paths.

A further examination yields additional information.  Consider our matching $M$ as a permuation $g: \{1,2,3,4,5\} \rightarrow \{A,B,C,D,E\}$.  And consider the function $f:  \{1,2,3,4,5\} \rightarrow \{A,B,C,D,E\}$ that maps the each path's start vertex to the fourth edge of the path.  So, for example, in the last example above, $f(1) = A, f(2) = A, f(4)=C, f(3) = D, f(5)=B$.  It is easy to see that each element in the target set has pre-image of cardinality 0, 1 or 2, so this will always be a 2-max function. This will be true in general. The algorithm will work if our matching $g$ has $g(i) \neq f(i)$ for each $ i = 1, ..., 5 $.  That is, the algorithm succeeds whenever the g-permutation is an f-derangement of the 2-max function $f$.  Therefore, by the analysis of the previous sections, we find it reasonable to conclude (with rough randomness assumptions on $f$) that a single iteration the algorithm will be successful with probability about $1/e$.  

We give one final observation.  When we use Konig's Theorem to get a 1-factor of $G$, that can easily be extended to get a decomposition of $G$ into 5 perfect matchings. Moreover, in Step 3 of the algorithm, we could just as easilty add the perfect matching edges to the end of each path, rather than the beginning.  Thus it is easy to see we have at least  5*2 = 10 options of how to run the algorithm, so should have confidence that at least one of these yields the desired decomposition.

Statements and Declarations:\\
The authors declare that no funds, grants or other support were received during the preparation of this manuscript.\\
The authors have no relevant financial or non-financial interests to disclose.\\
All authors contributed to the study conception and design.\\

Data Availability:\\
All data generated or analysed in this manuscript are included in this article.\\
\\
MSC codes 05A, 05C70


\begin{thebibliography} {10}
\bibitem{Botler} F. Botler, G. O. Mota and Y. Wakabayashi, Decompositions of triangle-free 5-regular graphs into paths of length five, {\em Discrete Math} 338 (2015) 1845-1855
\bibitem{Bouchet} A. Bouchet and J .L. Fouchet, Trois types de decompositions d'un graphe en cha\^{i}nes, Ann. Discrete Math 17 (1983) 131-141.
\bibitem{Graham}  R.L. H\"{a}ggkvist, Decompositions of complete bipartite graphs, in {\em Surveys in Combinatorics}, ed. J. Siemons, 1989 Cambridge University Press 115-146.
\bibitem{Jacobson} M. S. Jacobson, M. Truszczynski and Z. Tuza, Decompositions of regular bipartite graphs {\em Discrete Math}  89 (1991) 17-27.
\bibitem {Kotzig}  A. Kotzig, From the theory of finite regular graphs of degree three and four {\em \v{C}asopis P\v{e}st Mat} 82 (1957) 76-92.
\bibitem{Mont} R. Montgomery, A. Pokrovskiy and B. Sudakov,  A proof of Ringel's conjecture, {\em Geometric and Functional Analysis} 31 (2021) 663-720.
\bibitem{Ringel} G. Ringel, Problem 25. {\em Theory of Graphs and Its Applications}, Nakl. CSAV, Prague (1964)  162.
\bibitem{Wells}   D. Wells, {\em The Penguin Dictionary of Curious and Interesting Numbers} Middlesex, England, Penguin Books, p.27, 1986.
\bibitem{West} D. B. West, {\em Introduction to Graph Theory}, Prentice Hall, (2007)
\end{thebibliography}
\end{document}